# A Non-Existence Property of Pythagorean Triangles with a 3-D Application

*Konstantine Zelator*
*Department of Mathematics and*
*Computer Science*
*Rhode Island College*
*600 Mount Pleasant Avenue*
*Providence, RI 02908*
*USA*

E-Mail addresses:

kzelator@ric.edu
Konstantine_zelator@yahoo.com

## 1. Introduction

For students, teachers, and researchers alike, Pythagorean triangles hold a special fascination. Together with the wealth of interesting information about them, they have a history going back thousands of years. Oftentimes a result on this subject (and in mathematics in general) comes about after initially making a simple but critical observation. Such is the case with this work.

Generally speaking, it is not difficult to generate pairs of Pythagorean triangles (right triangles with integral side-lengths) such that one leg of one triangle is the hypotenuse of the other. In fact one can generate infinitely many such pairs. This can be easily done, and we do so in Section 5.

Now, consider the following question. Does there exist a pair of Pythagorean triangles such that the leg with the largest length in one triangle is the hypotenuse of the other; and in addition, the leg of smallest length in the first triangle, having the same length as a leg of the second? The answer to this question is in the negative; and it is proven in Theorem 1, in Section 6. This is the nonexistence property referred to in the title of this paper. The 3-D or three dimensional application also referred to in the title, is an application to Pythagorean boxes which are discussed in Section 7. The 3-D application is the content of Theorem 2, established in Section 8. Theorem 2 is an immediate consequence Proposition 2; and Theorem 1 is also an immediate consequence of Proposition 2 which itself follows from Proposition 1, in conjunction with Result 4. To establish Proposition 1, we make use of Results 1, 2, and 3. All four results; Results 1, 2, 3, and 4 are listed in Section 3. Results 1 and 2 sometimes are found as exercises in undergraduate number theory texts. On the other hand, results 3 and 4 are more advanced in nature, but well known in the literature. All four results are fully referenced at the end of this paper. The two propositions and their proofs are found in Section 4. The proofs in this work can be understood without any serious difficulties, by a good undergraduate student taking a number theory course. In Section 2 we present the familiar parametric formulas which generate the entire family of Pythagorean triples or triangles. Finally, in section 7, a geometric interpretation of Theorem 1 is offered. **Notation.** We use the standard notation *(a, b)* to denote the greatest common divisor of two integers *a and b.* When *(a, b) = 1,* the two integers are said to be *relatively prime*. The phrases, *"a* divides *b", "a* exactly divides *b", "b* is divisible by *a", "b* is a multiple of *a", "a* is a factor of *b";* all have the same meaning, namely that there exists an integer *k* such that *b = ka.* This fact can also be denoted by $a \mid b$. We also use the language of congruences in some places: $a \equiv b \pmod{n} \Leftrightarrow n \mid (a-b)$

### *2.* Pythagorean triples

We state the well-known parametric formulas describing the entire family of Pythagorean triples.

*The positive integers a, b, and c, satisfy* $a^2+b^2=c^2$, *if and only if* $a = \delta(m^2-n^2)$,
$b = \delta(2mn)$, $c = \delta(m^2+n^2)$;
*(or alternatively* $a = \delta(2mn)$ *and* $b = \delta(m^2-n^2)$*)* **(1)**
*for some positive integers* $\delta$, *m, n such that* $m>n$, $(m, n) = 1$,
*and* $m+n \equiv 1 \pmod 2$ *(which says that one of m, n is even, the other odd)*
*When* $\delta = 1$, *the Pythagorean triple (a, b ,c) is said to be* ***primitive****;*
*in which case* $(a, b) = (b, c) = (c, a) = 1$

For a wealth of historical information on the subject, the interested reader may refer to the book *"History of the Theory of Numbers, Vol. II"*; by L.E. Dickson (see **[1]**). Also, a great deal of material on Pythagorean triangles can be found in W. Sierpinski's "Elementary Theory of Numbers" (see **[2]**). And for a more textbook type of approach, one may refer to the book "*Elementary Number Theory and its applications*"; by Kenneth H. Rosen (see **[3]**).

## 3. Four Results from Number theory

The first two results below; Results 1 and 2 can be easily found in references **[2]** and **[3]**. Typically they are found in undergraduate number theory texts; Result 1 is usually listed as an exercise, while Result 2 as a theorem.

**Result 1** *Let $a, b, n$ be positive integers and suppose that $a^n|b^n$. Then, $a|b$. In particular, if $a^2|b^2$, then $a|b$.*

**Result 2** *Suppose that the positive integers a, b are relatively prime; $(a, b) = 1$. And that $ab = c^n$, for some natural number $n$. Then $a = a_1^n$ and $b = b_1^n$, for some positive integers $a_1$ and $b_1$ such that $(a_1, b_1) = 1$ and $a_1 b_1 = c$.*

The following two results below are more advanced in nature. Result 3 is historically attributed to P. Fermat. A detailed and carefully laid out proof of this result can be found in **[2].** The proof makes use of the method of infinite descent, which is based on the least element principle of the natural numbers.

**Result 3** *There exists no pair of natural numbers such that the sum of their squares is an integral square; and the difference of their squares is an integer square as well.*

The following result deals with the Diophantine equation $x^2 + 2y^2 = z^2$, which plays a central role in **[4]**, a paper published in the winter 1998 issue of *Mathematics and Computer Education*. This equation is a special case of the more general $x^2 + ky^2 = z^2$ , $k$ an integer with $k \geq 2$. A brief summary of this more general diophantine equation can be found in **[1]**. For a detailed, step-by-step analysis and derivation of its general solution, refer to **[5]**.

**Result 4** *The entire family of positive integer solutions of the Diophantine equation*

*$x^2 + 2y^2 = z^2$; can be described by the parametric formulas,*

$$x = \delta|k^2 - 2\lambda^2|;\ y = \delta(2k\lambda),\ z = \delta(k^2 + 2\lambda^2), \qquad \textbf{(2)}$$

*where $\delta, k, \lambda$ are positive integers such that $(k, \lambda) = 1$.*

*Note: It is not difficult to see, that these formulas reduce to the case $k \equiv 1(mod2)$.*

## 4. Two Propositions and Their Proofs

**Preposition 1** The diophantine equation $z^2 = x^4 + 4y^4$, has no solutions in positive integers x, y, z.

**Proof** Let, to the contrary, (x, y, z) be such a solution to $z^2 = x^4 + 4y^4$ **(3)**

We will show that this assumption leads to a contradiction.

Let d = (x, y). Then we must have,

$$\left\{\begin{matrix} x = dx_1, y = dy_1, with\ (x_1, y_1) = 1 \\ \text{for positive integers } x_1, y_1 \end{matrix}\right\} \quad \textbf{(4)}$$

From **(3)** and **(4)** it follows that $d^4|z^2$; and by Result 1, $d^2|z$.
Put $z = d^2 z_1$ for some natural number $z_1$. We obtain,

$$z_1^2 = x_1^4 + 4y_1^4 \quad \textbf{(5)}$$

Since $x_1$ and $y_1$ are relatively prime; they must either be both odd; or one odd, the other even. An argument modulo 8 easily shows they cannot be both odd. If there were both odd, then **(5)** shows that so would $z_1$. A well known fact from number theory is that the square of any odd integer; is congruent to 1 modulo 8. Thus we would have,
$z_1^2 \equiv x_1^4 \equiv y_1^4 \equiv 1\ (mod\ 8) \Longrightarrow (by\ (\mathbf{5}))\ 1 \equiv 1 + 4 \equiv 5\ (mod\ 8)$, impossibility. Thus $x_1$ and $y_1$ must have different parities. There are two possibilities or cases:

*Possibility (or Case)* 1: $x_1$ is odd, $y_1$ is even

*Possibility (or Case)* 2: $x_1$ is even, $y_1$ is odd

Below we show that the second possibility easily reduces to the first. Indeed, in the second possibility we have $x_1 = 2^{\propto} \bullet x_2, \propto \geq 1$, and with $x_2$ an odd integer. Substituting for $x_1$ in **(5)** yields,

$$z_1{}^2 = 2^2 \bullet (2^{4\propto-2} \bullet x_2^4 + y_1^4) \quad \textbf{(6)}$$

According to **(6)**, $z_1$ must be even; so that $z_1 = 2z_2$, for some positive integer $z_2$.

Equation **(6)** and a simple calculation produces,

$$z_2^2 = 2^2 \bullet (2^{\propto-1} \bullet x)^4 + y_1^4 \quad \textbf{(7)}$$

Note that $\propto \geq 2$; for if $\propto = 1$; **(7)** would be contradictory modulo 8, as we have seen above; by virtue of $z_2$, $x_2$, and $y_1$ all being odd. Thus, $\propto \geq 2$; $\propto -1 \geq 1$, and by putting $x_3 = 2^{\propto-1} \bullet x_2$, we see that **(7)** becomes $z_2^2 = 4 \bullet x_3^2 + y_1^4$; with $(x_3, y_1) = 1$ and $x_3$ even, $y_1$ odd; which is really under Possibility 1. This makes it clear that Possibility 2 above reduces to Possibility 1. Now, let us go back to equation **(5)** and consider Possibility 1. We have,

$$z_1^2 = (x_1^2)^2 + (2y_1^2)^2$$

Since $(x_1, y_1) = 1$ and $x_1$ is odd; it follows that $(x_1^2, 2y_1^2) = 1$. Therefore equation **(5)** describes a primitive Pythagorean triple with hypotenuse length $z_1$, odd leg length $x_1^2$ and even leg length $2y_1^2$. Accordingly, by formulas **(1)** we must have,

$$\left\{\begin{matrix} x_1^2 = m^2 - n^2, 2y_1^2 = 2mn, z = m^2 + n^2, \\ \text{for positive integers m, n such that } m > n, (m, n) = 1, and\ m + n \equiv 1 (mod\ 2) \end{matrix}\right\} \qquad \textbf{(8)}$$

The second equation in **(8)**, $y_1^2 = mn$, implies by Result 2 (with exponent $n = 2$) that each of $m, n$ must be an integer square: $m = m_1^2, n = n_1^2$. On the other hand, $x_1^2 = (m - n)(m + n)$, also implies that each of $m - n$ and $m + n$ must be an integer square; by Result 2 again. That is because $(m - n, m + n) = 1$, which easily follows from the conditions $(m, n) = 1$ and $m + n \equiv 1 (mod\ 2)$ in **(8)** (This is a standard exercise in a number theory course). Putting everything together, we have $m_1^2 - n_1^2 = k^2$ and $m_1^2 + n_1^2 = \ell^2$, for positive integers $m_1, n_1, \ell, k$; which is in violation of Result 3. The proof is complete ■.

**Proposition 2** The four-variable Diophantine system $\left\{\begin{matrix} z^2 = w^2 + y^2 \\ w^2 = y^2 + x^2 \end{matrix}\right\}$, **(9)** has no solutions in positive integers $x, y, z, w$.

**Proof** If to the contrary, such a solution $(x, y, z, w)$ existed. Then **(9)** implies, $z^2 = x^2 + 2y^2$ **(10)**

By Result 4, equation **(10)** implies that, $x = \delta|k^2 - 2\lambda^2|, y = \delta(2k\lambda), z = \delta(k^2 + 2\lambda^2)$, **(11)**

for some positive integers $\delta, k, \lambda$ such that $(k, \lambda) = 1$. From **(11)** and the second equation in **(9)** it follows that $w^2 = \delta^2\ [(2k\lambda)^2 + (k^2 - 2\lambda^2)^2] \Leftrightarrow w^2 = \delta^2 \bullet (k^4 + 4\lambda^4)$; which turn implies (by Result 1) that $\delta | w$; put $w = \delta \bullet w_1$, for some natural number $w_1$.

The last equation above then yields, $w_1^2 = k^4 + 4\lambda^4$ **(12)**

Equation **(12)** says that $(k, \lambda, z)$ is a positive integer solution of the diophantine equation $z^2 = x^4 + 4y^4$, contradicting Proposition 1. End of proof ■.

## 5. Pythagorean Triangles with a Common Side

Picture two Pythagorean triangles sharing a side. We leave it to the interested reader to explore the case when the common side is a leg. In this section, we focus on the case wherein the hypotenuse of one Pythagorean triangle is a leg of the other. If we pick two Pythagorean triples $(a_1, b_1, c_1)$ and $(a_2, b_2, c_2)$ at random. Then by **(1)**,

Pythagorean triangle 1: $a_1 = \delta(m^2 - n^2), b_1 = \delta(2mn), c_1 = \delta(m^2 + n^2)$
Pythagorean triangle 2: $a_2 = d(M^2 - N^2),\ b_2 = d(2MN),\ c_2 = d(M^2 + N^2)$;
where $m, n$ satisfy the conditions in (1); and likewise for $M$ and $N$.

In order for the hypotenuse of Pythagorean triangle 2 to be a leg of Pythagorean triangle 1 we must have either, $\delta\ (2mn) = d(M^2 + N^2)$ **(13a)**

Or alternatively $\delta\ (m^2 - n^2) = d(M^2 + N^2)$ **(13b)**

Indeed, it is easy to generate entire families of such pairs of Pythagorean triangles (the interested reader may want to explore this further). For example to make (**13a**) satisfied take,

Family 1: $d = 4,\ n = 2, \delta = 1,\ m = M^2 + N^2$
Family 2: $\delta = K \bullet (\mathrm{M}^2 + \mathrm{N}^2),\ d = K(2mn),\ K$ a positive integer

To make (**13b**) satisfied take,

Family 3: $\delta = K(M^2 + N^2), d = K \bullet (m^2 - n^2),\ K$ is a natural number.
Family 4: $\delta = 1,\ d$ odd, $(d, M^2 + N^2) = 1,\ d > M^2 + N^2$;
And with $m = \frac{d + M^2 + N^2}{2},\ n = \frac{d - (M^2 + N^2)}{2}$

## 6. Theorem 1 and its Proof

**Theorem 1** There exists no pair of Pythagorean triangles such that the leg of largest length in the first triangle is the hypotenuse of the second triangle; and in addition, with the leg of smallest length in the first triangle; having the same length as a leg in the second triangle.

**Proof** If two such Pythagorean triples $(a, b, c)$ and $(a, d, b)$ existed; we would have

$\begin{Bmatrix} c^2 = b^2 + a^2 \\ b^2 = d^2 + a^2 \end{Bmatrix}$, which says that the 4-tuple $(a, d, b, c)$ is a positive integer solution to the system $\begin{Bmatrix} z^2 = w^2 + y^2 \\ w^2 = y^2 + x^2 \end{Bmatrix}$, contrary to Proposition 2. ■

## 7. Pythagorean boxes

A *Pythagorean box* is a rectangular parallelepiped whose edges have integer lengths; and with the two inner diagonals also having integral length. Such a solid has three pairs of congruent faces, with each face being a rectangle. The twelve edges can be divided into three groups, with each group containing four sides of common length. Let $x, y$, and $z$ be these three integral lengths. And let $t$ be the length of each of the two inner diagonals. Then, as it can be easily seen from the obvious geometry, the four positive integers $x, y, z$, and $t$ must satisfy the equation, $t^2 = x^2 + y^2 + z^2$ **(14)**

The general solution of the Diophantine equation **(14)** is well known and it can be found (together with

Its derivation) in **[2].** In fact the entire family of solutions of 14; in positive integers $x, y, z, t$ ; can be parametrically described by the formulas, $x = 2v, y = 2\ell, z = \frac{\ell^2+v^2-n^2}{n}, t = \frac{\ell^2+v^2+n^2}{n}$, **(14a)** where $v, \ell, n$ are positive integers such that $n$ is a divisor of $\ell^2 + v^2$; and $n < \sqrt{\ell^2 + v^2}$. If we take $n = 1$ in (14a), we obtain the subfamily of solutions,

$$x = 2v, \ \ y = 2\ell, \ \ z = \ell^2 + v^2 - 1, \ \ t = \ell^2 + v^2 + 1 \qquad \textbf{(14b)}$$

We can use **(14b)** to generate an infinite subset of solutions which produce Pythagorean boxes which have a face diagonal of integral length (and thus four congruent face diagonals having the same integer length). To do so, all we have to do is require that $x^2 + y^2 = d^2$; put $d = 2D$, which gives

$(2v)^2 + (2\ell)^2 = 4D^2$; $v^2 + \ell^2 = D^2$. And by the formulas in **(1)**,

$$v = \delta(m^2 - k^2), \ \ \ell = \delta(2mk), d = 2\delta(m^2 + k^2), \qquad \textbf{(14c)}$$

for positive integers $m, k$ such that $m > k$, $(m, k) = 1$, and $m + k \equiv 1\ (mod 2)$
Formulas **(14b)** and **(14c)** describe an infinite set of Pythagorean boxes, with the two congruent faces (with sides of lengths $x$ and $y$) each having two diagonals of length $d$. On the other hand, if we wish to generate Pythagorean boxes in which two adjacent edges have the same length; say $x = y$; in fact we can find all such Pythagorean boxes from **(2)**, since according to **(14)** we would have $t^2 = z^2 + 2x^2$ (In **(2),** replace $z$ by $t, x$ by $z$ and $y$ by $x$). The natural question to ask at this stage is this: Are there Pythagorean boxes of whose two congruent faces are squares; while another pair of congruent faces have diagonals of integer length? Theorem 2 below provides the answer in the negative.

### 8. Theorem 2 and its Proof

**Theorem 2** There exists no Pythagorean box with a pair of congruent (or opposite) faces being squares. And with the four diagonals of equal length in another pair of opposite faces, also having integral length.

**Proof** If such a Pythagorean box existed**;** represented by $(r, s, p, q)$. We would have $q^2 = r^2 + s^2 + p^2, r = s$, and $d^2 = r^2 + p^2$; for positive integers $r, s, p, q$. This then easily implies, $q^2 = d^2 + r^2$ and $d^2 = r^2 + p^2$; showing that $(p, r, q, d)$ is a positive integer solution to the system $\begin{Bmatrix} z^2 = w^2 + y^2 \\ w^2 = y^2 + x^2 \end{Bmatrix}$; contrary to Proposition 2. ■